\documentclass[12pt]{article}

\usepackage{amsmath}
\usepackage{amsfonts}
\usepackage{amsthm}
\usepackage{graphicx}
\usepackage{hyperref}
\usepackage{enumerate}
\usepackage[parfill]{parskip}

\numberwithin{equation}{section} 

\newtheorem{thm}{Theorem}[section]

\newtheorem{conj}[thm]{Conjecture}

\newcommand{\N}{\mathbb{N}}         

\begin{document}

\title{Explicit expressions for the moments of the size of an $(n,dn-1)$-core partition with distinct parts}

\author{
Anthony Zaleski\thanks{Department of Mathematics, Rutgers University (New Brunswick), 110 Frelinghuysen Road, Piscataway, NJ 08854-8019, USA.}
}

\maketitle
\begin{abstract}
In a previous paper (arXiv:1608.02262), we used computer-assisted methods to find explicit expressions for the moments of the size of a uniform random $(n,n+1)$-core partition with distinct parts.  In particular, we conjectured that the distribution is asymptotically normal.  However, our analysis hinged on a characterization of $(n,n+1)$-core partitions given by Straub, which is not readily generalized to other families of simultaneous core partitions.

In another paper (arXiv:1611.05775) with Doron Zeilberger, we made use of the characterization in terms of posets to analyze $(2n+1,2n+3)$-core partitions with distinct parts; here, the distribution was found \emph{not} to be asymptotically normal.  Our method involved finding recursive structure in the relevant sequence of posets.  We remarked that this method is applicable to other families of core partitions, provided that one can understand the corresponding posets. 

Here, we use the poset method (and, as before, a computer) to analyze $(n,dn-1)$-core partitions with distinct parts, where $d$ is a natural number.  (This problem was introduced by Straub in arXiv:1601.07161.)  We exhibit formulas for the moments of the size, as functions of $d$ with $n$ fixed, and vice versa.  We conjecture that the distribution is asymptotically normal as $n$ approaches infinity.  Finally, we find formulas for the first few moments, as functions of \emph{both} $n$ and $d$.
\end{abstract}


\\
\\
\section{Introduction}\label{sec:intro}

\begin{figure}[ht]
\begin{center}
\includegraphics[width=.5\textwidth]{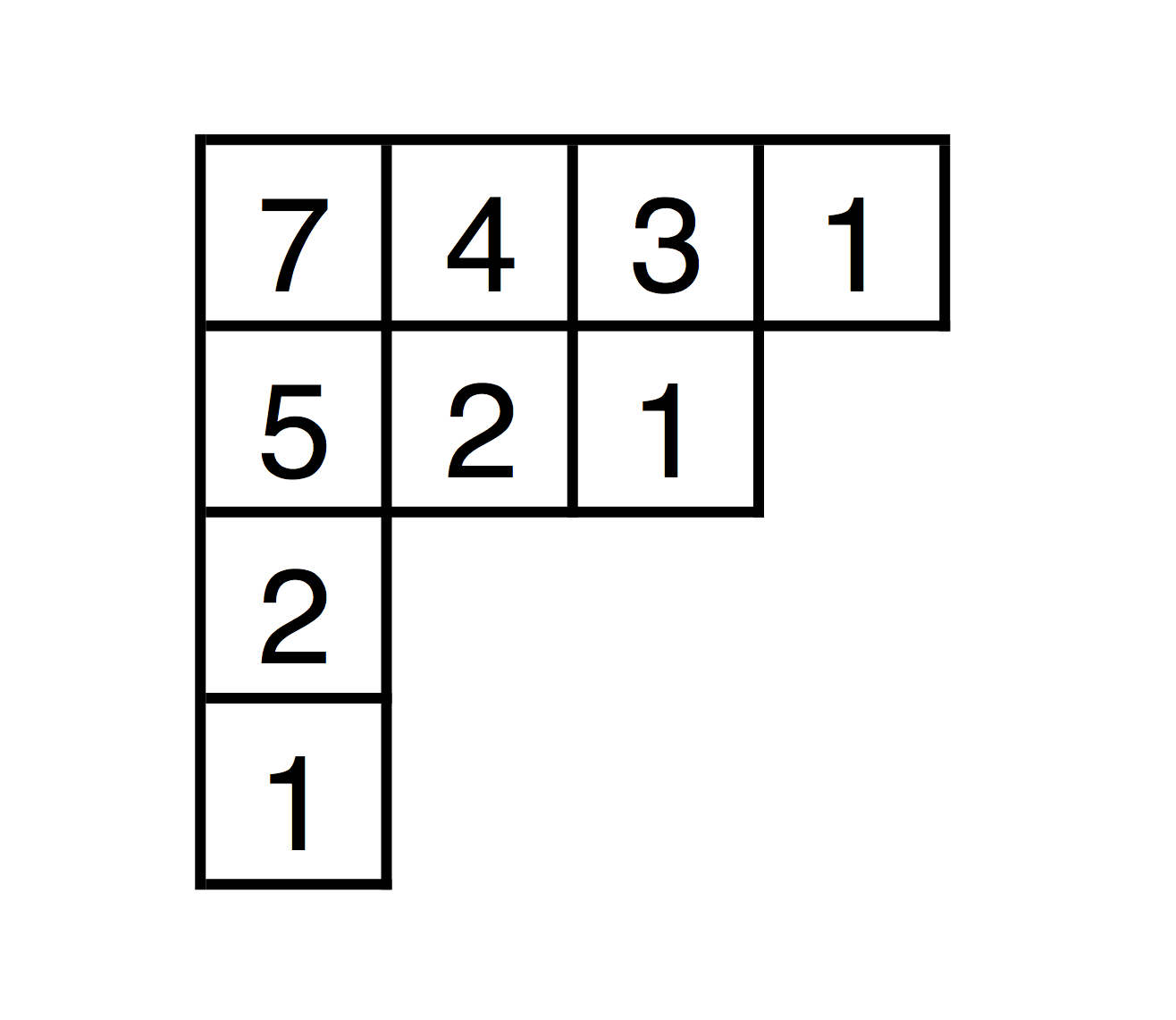}
\caption{Young diagram of the partition $9=4+3+1+1,$ showing the hook lengths of each box.}
\label{fig:young}
\end{center}
\end{figure}

\subsection{$(s,t)$-core partitions}
Recall that the \emph{hook length} of a box in the Young diagram of a partition is the number of boxes to the right (the arm) plus the number of boxes below it (the leg) plus one (the head).  (We use the English convention for Young diagrams; see Figure \ref{fig:young}.)   A partition is an $s$-\emph{core} if its Young diagram avoids hook length $s$ and an $(s,t)$-\emph{core} if it avoids hook lengths $s$ and $t$ [AHJ].  For example, the partition $9=4+3+1+1$ in Figure 1 is a $(6,8)$-core but not a $(6,7)$-core.

The number of $(s,t)$-core partitions is finite iff $s$ and $t$ are coprime, which we shall assume from now on [AHJ].  Let $X_{s,t}$ be the random variable ``size of an $(s,t)$-core partition," where the sample space is the set of all $(s,t)$-core partitions, equipped with the uniform distribution.  In [EZ], with the help of Maple, Zeilberger derived explicit  expressions (as polynomials in $s$ and $t$) for the expectation, variance, and numerous higher moments of $X_{s,t}$.  The original paper noted that ``From the `religious-fanatical' viewpoint of the current `mainstream' mathematician, they are `just' conjectures, but nevertheless, they are \textbf{absolutely certain} (well, at least as absolutely certain as most proved theorems)," and a donation to the OEIS was offered for the theory to make the results rigorous.  Later, it was found that such theory did exist and the results are entirely rigorous; see the updates at the paper's site.  

Zeilberger also computed some standardized central moments of $X_{s,t}$ and the limit of these expressions as $s,t\to \infty$ with $s-t$ fixed.  From this he conjectured the limiting distribution.  Perhaps surprisingly, it is abnormal.
\subsection{$(n,n+1)$-core partitions with distinct parts}
Things are not as easy if look at the random variable  ``size of an $(s,t)$-core partition \emph{with distinct parts}.''  In this case, there does not seem to a single formula for the moments in terms of $s$ and $t$.  However, we can consider certain indexed families of core partitions and look for explicit formulas for the moments.

For example, in [S], Straub remarks in a lemma that a partition with distinct parts is an $(n,n+1)$-core iff it has perimeter (largest hook length) less than $n$. Using this, it is easy to show that the number of such partitions is the Fibonacci number.  Further, in [Z], we used Straub's lemma and some experimentation with Maple to derive explicit expressions for the moments, as rational functions of Fibonacci numbers.  Our results give strong evidence that the limiting distribution as $n\to\infty$ \emph{is} normal in this case.
\subsection{The poset characterization}
Unfortunately, we do not see how to generalize the characterization above to other indexed families of partitions.  However, there is another characterization of $(s,t)$-core partitions due to Anderson [A].

Given coprime $s$ and $t$, define the poset 
\begin{equation}\label{eq:Pst}
P_{s,t} := \N \setminus (s\N + t\N),
\end{equation}
where $\N=\{ 0,1,2,3, \dots, \}$, where the partial-order relation $c \, \leq_P \, d$ 
holds whenever $d-c$ can be expressed as $\alpha s + \beta t$ for some $\alpha,\beta \in \N$.  (It can be shown that $P_{s,t}$ is finite with largest element $st-s-t$.)

Anderson defined a bijection between $(s,t)$-core partitions and \emph{order ideals} of $P_{s,t}$.  (Under our convention, a subset  $I$ of a poset $P$ is an order ideal iff $x \in I, y\leq_P x$ $\implies $ $y \in I$.)  Further, the hook lengths of the boxes in the leftmost column in the partition's Young diagram correspond to the labels of the order ideal.  

Now, a partition has distinct parts iff there are no consecutive hook lengths in the leftmost column. Thus:
\begin{itemize}
\item[] \emph{$(s,t)$-core partitions with distinct parts are bijective to order ideals of $P_{s,t}$ with no consecutive labels.}
\end{itemize}
\subsection{$(2n+1,2n+3)$-core partitions with distinct parts}
So, given a family of core partitions, we have a corresponding sequence of posets which we can hopefully understand.  In [ZZ], we plotted some of the posets corresponding to $(2n+1,2n+3)$-core partitions with distinct parts and were able to see a recursive structure.  This led to recursions for the generating functions of the partitions.  We were again able to derive explicit expressions for the moments in terms of $n$ and show that the distribution of size is not asymptotically normal as $n\to\infty$.
\section{$(n,dn-1)$-core partitions with distinct parts}
In [S], Straub generalizes the problem in Section 1.2 by considering $(n,dn-1)$-core partitions with distinct parts, where $n$ and $d$ are natural numbers. He proves in Theorem 4.1 that the number of such partitions, call it $N_d(n)$, satisfies a generalized Fibonacci recurrence:
\begin{align}
\begin{split} \label{eq:Nds}
&N_d(1)=1, \,\, N_d(2)=d, 
\\
&N_d(n)=N_d(n-1)+dN_d(n-2).
\end{split}
\end{align}
Of course, this reduces to the usual Fibonacci numbers when $d=1$.  Note that we can view $N_d(n)$ as a sequence of polynomials in $d$.

Here, we shall use the poset characterization in Section 1.3 to easily recover Straub's result and discover new conjectures about the distribution of the sizes of the partitions.

\subsection{Understanding the posets}
By Section $1.3$, we know that $(n,dn-1)$-core partitions with distinct parts are bijective with order ideals of $P_{n,dn-1}$ containing no consecutive labels.  We can use the procedure \verb+PW+ in the Maple package \verb+Armin+ accompanying [ZZ] to plot $P_{n,dn-1}$ for various $n$ and $d$.

\begin{figure}[ht]
\begin{center}
\includegraphics[width=.6\textwidth]{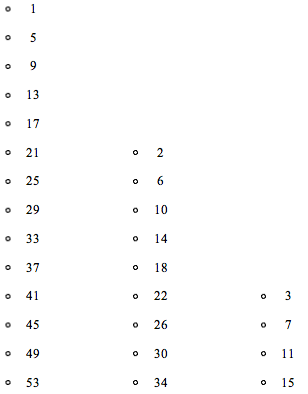}
\caption{The poset $P_{4,19}=P_{4,4\cdot5-1}$.}
\label{fig:pw}
\end{center}
\end{figure}

For example, Figure \ref{fig:pw} depicts the poset $P_{4,19}$, i.e., the $n=4, d=5$ case.  (When plotting $P_{s,t}$, we use the convention to increment $s$ in the $\downarrow$ direction and $t$ in the $\leftarrow$ direction.  Thus, the largest label, $st-s-t$, is in the lower left corner.)  It is easy to show that this general trapezoidal shape persists for arbitrary values of $n$ and $d$.  We can also see $P_{n,dn-1}$ as a colonnade of $n-1$ vertical pillars with heights $d(n-1)-1,d(n-2)-1,\dots,d-1$.  Further, the tops of the pillars have labels $1,2,\dots,n-1$.

\subsection{Characterizing the order ideals}
Next,  we recover the recursion \eqref{eq:Nds} by enumerating the order ideals of $P_{n,dn-1}$ with no consecutive labels.  

Referring to the $n=4, d=5$ example, let $I$ be an order ideal of $P_{4,19}$.  Let $I_k$ be the part of $I$ contained in the $k$th pillar.  
\begin{itemize}
\item
If $I_1=\emptyset$, then $I$ is isomorphic to an order ideal of $P_{n-1,d(n-1)-1}=P_{3,5\cdot 3-1}$.  
\item
Otherwise, let $x$ be the largest member of $I_1$.  Then $x \in \{1,1+n,\dots,1+(d-1)n\}=\{1,5,9,13,17\}$. For if $21 \in I_1$, then $1,2\in I$, contradicting the assumption that $I$ contains no consecutive labels.  So there are $d=5$ choices for $I_1$.  Further, $I_2$ must be empty; otherwise, again, we have $1,2 \in I$.  Thus the remainder of $I$ is isomorphic to an order ideal of $P_{n-2,d(n-2)-1}=P_{2,5\cdot2-1}$.
\end{itemize}

To summarize, if $I$ is an order ideal of $P_{n,dn-1}$ with no consecutive labels, then either $I$ is isomorphic to an order ideal of $P_{n-1,d(n-1)-1}$, or $I_1$ has $d$ options and the rest of $I$ is isomorphic to an order ideal of $P_{n-2,d(n-2)-1}$.  This proves \eqref{eq:Nds}.

From the above observation, we have the following characterization:

\begin{itemize}
\item[]Any order ideal of  $P_{n,dn-1}$ with no consecutive labels is of the form $I=I_1\cup\cdots\cup I_{n-1}$, where
\begin{itemize}
\item
 Each $I_k$ is either empty or of the form $\{k,k+n,\cdots,k+i_kn\}$, where $i_k \leq d-1$ if $k<n-1$, and $i_{n-1}<d-1$.
 \item
 If $I_k$ is nonempty, then $I_{k+1}$ is empty.
 \end{itemize}
 \end{itemize}
 
 In short: To make an order ideal, we hang strings of beads from the tops of the pillars in such a way the strings are not too long and adjacent pillars are not both decorated.
 \subsection{Computing the generating function}
 Our ultimate goal is to investigate the distribution of the size of $(n,dn-1)$-core partitions with distinct parts.  To this end, we define the generating function
 \begin{equation}
 G_{d, n}(q):=\sum_p q^{|p|},
 \end{equation}
 where $p$ ranges over $(n,dn-1)$-core partitions with distinct parts, and $|p|$ denotes the size of the partition $p$, i.e., the sum of its parts.  We shall give an efficient scheme for computing $G_{d,n}(q)$ for fixed $d$ and $n$.
 
Proceeding as in [ZZ], we first compute the auxiliary generating function
 \begin{equation} \label{eq:Fdn}
 F_{d, n}(q,t):=\sum_I q^{w(I)}t^{|I|},
 \end{equation} 
 where $I$ ranges over all order ideals of $P_{n,dn-1}$ with no consecutive labels;  $w(I)$ is the sum of the labels in $I$; and $|I|$ is the number of labels in $I$.  Then, as explained in [ZZ], we can obtain $G_{d, n}(q)$ by replacing occurrences of $t^k$ in $F_{d, n}(q,t)$ with $q^{-k(k-1)/2}$.
 
 To compute $F_{d, n}(q,t)$, we use the reasoning of the previous section, but this time we keep track of the weight of the order ideal.  
 
 First, we introduce yet another auxiliary generating function.  For $1\leq k \leq n-1$, let $P^k_{n,dn-1}$ be the sub-poset  of $P_{n,dn-1}$  obtained by chopping off everything to the left of the $k$th column (note $P^1_{n,dn-1}=P_{n,dn-1}$).  Define  $F^k_{d, n}(q,t)$ as in \eqref{eq:Fdn}, except $I$ ranges over order ideals of $P^k_{n,dn-1}$ with no consecutive labels (note $F^1_{d, n}(q,t)= F_{d, n}(q,t)$).
 
 By the reasoning of the previous section, the first column of an order ideal of $P^k_{n,dn-1}$ is either empty or of the form $\{k,k+n,\cdots,k+in\}$, where $0\leq i\leq d-1$).  Since the latter set has weight 
 $$
 q^{\sum_{j=0}^i (k+jn)}t^{i+1} = q^{(i+1)(in/2+k)}t^{i+1},
 $$
 we have the recursion
 \begin{align}
 \begin{split}  \label{eq:Fkdn}
 &F^k_{d, n}(q,t)=F^{k+1}_{d, n}(q,t)+\left(\sum_{i=0}^{d-1} q^{(i+1)(in/2+k)}t^{i+1} \right) F^{k+2}_{d, n}(q,t) \text{ for } 1\leq k \leq n-2;
 \\
 &F^{n-1}_{d, n}(q,t)=\sum_{i=0}^{d-2} q^{(i+1)(in/2+k)}t^{i+1};
 \\
 &F^{n}_{d, n}(q,t):=1.
 \end{split}
 \end{align}
 
 Note that this is a recursion in the auxiliary index $k$, not in $n$ and $d$.
 
 Given $n$ and $d$, we can use \eqref{eq:Fkdn} to find $F^1_{d, n}(q,t)= F_{d, n}(q,t)$.  Finally, we make the substitution $t^k \to q^{-k(k-1)/2}$ to find $G_{d,n}(q)$.  All of this is done in the procedure \verb+Gdn+ in the Maple package.
 \section{Distribution of the size}
 Given fixed $n$ and $d$, we can pick a uniform random $(n,dn-1)$-core partition with distinct parts, and consider its its size, call it $X_{n,d}$.  Then $X_{d,n}$ is a random variable, so it makes sense to inquire about its distribution.  Since $G_{d,n}$ is the  generating function for $X_{d,n}$, we can easily compute as many moments of the distribution as we please, for fixed $n$ and $d$.
 
Using this information, we can investigate how the moments behave as functions of $n$ and $d$.  We consider two cases: $n$ is variable and $d$ is fixed, and vice versa.  In each case, we consider the behavior of $X_{d,n}$ as the variable tends to infinity; in particular, we address the question of asymptotic normality. Finally, we derive formulas for the first few moments as functions of \emph{both} $n$ and $d$.
  \subsection{Varying $d$ and fixing $n$}
  First, we introduce some terminology. Given a natural number $k$, let us denote 
  $$
 m_k(d,n):=\left[ \left(q\frac{d}{dq}\right)^kG_{d,n}(q)\right]_{q=1}
  $$
  the $k$th ``pre-moment'' of $X_{d,n}$.  Thus,
  $$
  M_k(d,n):=\frac{m_k}{G_{d,n}(1)}=\frac{m_k}{N_d(n)}=\mathbb{E}[X_{d,n}^k]
  $$
  is the $k$th (straight) moment of $X_{d,n}$.  For example, the mean is $\mu_{d,n}=M_1$. 
  
  We denote the $k$th \emph{central moment} by
  $$
  M_k^c(d,n):=\mathbb{E}[(X_{d,n}-\mu)^k].
  $$
  For example, the variance is $\sigma_{d,n}^2:=M_2^c$.  
  
  Finally, the $k$th \emph{standardized moment} is
  $$
  M_k^s(d,n):=\frac{M_k^c}{\sigma^k}.
  $$
  Note that the central, straight, and standardized moments can easily be computed from the pre-moments.
  
  Now, for numeric values of $d$ and $n$, we can easily compute all the quantities above.  Analyzing moment data for many values of $d$ and $n$ confirms the following:
  
  \begin{conj}
  For each $n$, the $k$th pre-moment $m_k(d,n)$ of $X_{d,n}$ is a polynomial in $d$.  Further, the degree of this polynomial is $2k+\lfloor n/2 \rfloor$.
  \end{conj}


For example, our experimental evidence indicates that 

\begin{align}
\begin{split} \label{eq:lim}
\{\lim_{d\to\infty} M_k^s(d,3)\}_{k=3}^\infty &=
2/7\,\sqrt {5},{\frac {15}{7}},{\frac {100}{77}}\,\sqrt {5},{\frac {
6625}{1001}},{\frac {750}{143}}\,\sqrt {5},\dots
\\
& \approx .641, 2.14, 2.91, 6.62, 11.7,\dots
\\
\{\lim_{d\to\infty} M_k^s(d,4)\}_{k=3}^\infty &\approx
.162, 2.08, 1.19, 6.20, 7.05,\dots
\\
\{\lim_{d\to\infty} M_k^s(d,5)\}_{k=3}^\infty &\approx
.237, 2.22, 1.76, 7.43, 10.8,\dots
\\
\{\lim_{d\to\infty} M_k^s(d,6)\}_{k=3}^\infty &\approx
.052, 2.36, .671, 7.80, 5.15,\dots
\\
&\dots
\\
\{\lim_{d\to\infty} M_k^s(d,10)\}_{k=3}^\infty & \approx
-0.001, 2.62, .130, 10.1, 2.17,\dots.
\end{split}
\end{align}

Recall that the standard normal distribution has moments $0,1,0,3,0,15,\dots$.  The sequences above seem to approach this as $n\to\infty$, leading us to the following:
 \begin{conj}
 For each fixed $n$, the distribution of $X_{d,n}$ is \emph{not} asymptotically normal; in fact, $(X_{d,n}-\mu_{d,n})/\sigma_{d,n}$ tends to some abnormal distribution $X_n$ as $d\to\infty$.  However, $X_n$ \emph{is} asymptotically normal; that is, $(X_n-\mu)/\sigma$ tends to the standard normal distribution as $n\to\infty$.
 \end{conj}
 \subsection{Varying $n$ and fixing $d$}
 Next, we fix $d$, and look at $X_{d,n}$ as a sequence of random variables indexed by $n$.  The $d=1$ case was already addressed in [ZZ], where we found that the pre-moments are given by polynomials in $n$ and the Fibonacci numbers.  In light of \eqref{eq:Nds}, we might expect the same to be true for arbitrary $d$, except we use the generalized Fibonacci numbers, $N_d(n)$:

 \begin{conj}
  For each $d$, the $k$th pre-moment $m_k(d,n)$ of $X_{d,n}$ is of the form $a(n)N_d(n)+b(n)N_d(n+1)$, where $a$ and $b$ are polynomials in $n$.
  \end{conj}
  
Again, experimental evidence verifies this claim.  The one anomalous case seems to be $d=2$, for which $N_d(n)=2^{n-1}$.  In this case, our methods do not yield nice formulas for the moments.  

Upon computing the limits standardized moments, we do get the familiar sequence $0,1,0,3,0,15,\dots$ in this case, leading to the following:

\begin{conj} 
For fixed $d$, the distribution of $X_{d,n}$ is asymptotically normal.  That is, $(X_{d,n}-\mu_{d,n})/\sigma_{d,n}$ approaches the standard normal distribution as $n\to \infty$.
\end{conj}
 \subsection{Varying $n$ and $d$}
 Finally, it is possible (but computationally taxing) to obtain a single formula for the $k$th moment as a function of \emph{both} $k$ and $d$.  In the previous section, we \emph{fixed} $d$ and $k$, we and used \verb+GuessRecPol+ to fit the sequence $\{m_k(d,n)\}_{n=2}^\infty$ to the ansatz $a(n)N_d(n)+b(n)N_d(n+1)$, where $a(n)$ and $b(n)$ are polynomials.  However, using the methods Section 3.1, we can fix \emph{only} $k$ and look at $\{m_k(d,n)\}_{n=2}^\infty$ as a sequence of polynomials in $d$.  Further,  due to Maple's ability to handle linear systems with symbolic coefficients, \verb+GuessRecPol+ can still be applied to this sequence.  Of course, now $a(n)$ and $b(n)$ will be polynomials in $n$ whose coefficients are rational functions in $d$:

 \begin{conj}
  The $k$th pre-moment $m_k(d,n)$ of $X_{d,n}$ is of the form $A(n,d) N_d(n) + B(n,d)N_d(n+1)$, where $A$ and $B$ are degree $2k$ polynomials in $n$ whose coefficients are rational functions in $d$.
  \end{conj}
  
 Due to the amount of data needed to fit the $k$th moment to the ansatz, it takes a few minutes even to generate the formula for the 3rd moment.  
 \subsection{Sample formulas for the moments}
 Here, we present a small taste of the conjectures yielded by our Maple package.  In the first two conjectures to follow, we could easily have presented formulas for many more moments, but we omit them to save space.  See the next section for instructions to replicate these results and many more for yourself.
 
 In general, we conjectured that $M_k(d,n)$ is a rational function in $d$ for $n$ fixed.  However, for $n=3$ the straight moments seem to be polynomials:
 
 \begin{conj}
 The expectation of $X_{d,3}$ is $d^2/3+d/4-1/12$, and the variance is $4d^4/45+d^3/12-d^2/144+d/24+31/720$. 
\end{conj}
 
 Here is an example where we fix $d$.  
 
 \begin{conj}
 The expectation of $X_{3,n}$ is
 $$
 {\frac {25}{39}}\,{n}^{2}-{\frac {479}{507}}\,n+{\frac {406}{507}}+\frac{N_3(n+1)}{N_3(n)}\, \left( -\frac{1}{39}{n}^{2}+{\frac {29}{169}}\,n-{\frac {158}{
507}} \right).
 $$
 \end{conj}
 
 Finally, here is the expectation once and for all, in terms of both $n$ and $d$:
 
 \begin{conj}
 The expectation of $X_{d,n}$ is
 \begin{align*}
 &{\frac { \left( 5\,{d}^{3}+7\,{d}^{2}+d-1 \right) {n}^{2}}{24(4\,d+
1)}}-{\frac { \left( 8\,{d}^{5}+21\,{d}^{4}+7\,{d}^{3}-{d}^{2}+3
\,d-2 \right) n}{24(16\,{d}^{3}-24\,{d}^{2}-15\,d-2)}}
\\
&+{\frac {17\,{
d}^{4}+13\,{d}^{3}-9\,{d}^{2}-7\,d-2}{12(16\,{d}^{3}-24\,{d}^{2}-15\,d-2)}
}+\frac{N_d(n+1)}{N_d(n)}\, 
\\
&\cdot \left( -{\frac { \left( {d}^{2}-1 \right) {n}^{
2}}{24(4\,d+1)}}-{\frac { \left( 2\,{d}^{4}-9\,{d}^{3}-16\,{d}^{2}-3
\,d+2 \right) n}{8(16\,{d}^{3}-24\,{d}^{2}-15\,d-2)}}-{\frac {{d}^{
4}+20\,{d}^{3}+9\,{d}^{2}-20\,d-10}{ 12\left( d-2 \right)  \left( 4\,d+1
 \right) ^{2}}} \right).
\end{align*}
 \end{conj}
 
 Note that this formula is singular at $d=2$, explaining anomaly mentioned earlier.  However, we can still make sense of the $d=2$ case by first plugging in a numeric value of $n$, (so that the $N_d$'s become polynomials in $d$), then taking the limit as $d\to 2$.  So this formula effectively works for \emph{all} $n$ and $d$.
  
\section{Using the Maple package}
The Maple package \verb+core2.txt+ accompanying this paper may be found at the following URL:
 \\ \url{http://www.math.rutgers.edu/~az202/Z}.  
 
 To use the Maple package, place \verb+core2.txt+ in the working directory and execute \verb+read(`core.txt`);+.  
 
 To see the main procedures, execute \verb+Help();+.  For help on a specific procedure, use \verb+Help(<procedure name>);+.  
 
 Here are some things to try:
 \begin{itemize}
 \item \verb+Gdn(q,3,7);+ gives the generating function (according to size) of $(3,3\cdot 7-1)$-core partitions with distinct parts.
 \item \verb+MdnK(d,3,1);+ and \verb+McdnK(d,3,2);+ reproduce Conjecture 3.6.
 \item \verb+MdnK(3,n,1);+ reproduces Conjecture 3.7.
 \item \verb+MdnK(d,n,1);+ reproduces Conjecture 3.8.
 \item \verb+map(p->limit(p,d=infinity),MsdnK(d,3,7));+ reproduces the first equation in \eqref{eq:lim}.
 \end{itemize}
\section*{Acknowledgements} 
The author thanks Dr. Doron Zeilberger for introducing this project to him and  guiding his research in the right direction. The author also thanks Dr. Neil Sloane for proofreading the draft.
\section{References}
\begin{itemize}
\item[{[A]}]
J. Anderson, Partitions which are simultaneously $t_1$- and $t_2$-core, Discrete Math. 51 (2014) 205-220.
\item[{[AHJ]}]
D. Armstrong, C. R. H. Hanusa, B. C. Jones, Results and conjectures
on simultaneous core partitions, European J. Combin. 41 (2014) 205-
220.
\item[{[EZ]}]
S. B. Ekhad, D. Zeilberger, Explicit expressions for the variance and higher moments of the size of a simultaneous core partition and its limiting distribution, The Personal Journal of Shalosh B. Ekhad and Doron Zeilberger, posted Aug. 30, 2015,
\\
\url{http://www.math.rutgers.edu/~zeilberg/mamarim/mamarimhtml/stcore.html}.
\item[{[S]}]
A. Straub, Core partitions into distinct parts and an analog of Euler's theorem, European J. Combin. 57 (2016) 40-49.
\item[{[Z]}]
A. Zaleski, Explicit expressions for the moments of the size of an $(s,s+1)$-core partition with distinct parts, Adv. Appl. Math. 84 (2017) 1-7.
\item[{[ZZ]}]
A. Zaleski, D. Zeilberger, Explicit (polynomial!) expressions for the expectation, variance, and higher moments of the size of a $(2n+1,2n+3)$-core partition with distinct parts, preprint, arXiv:1611.05775, 2017.
\end{itemize}
\end{document}